\documentclass[12pt]{article}
\usepackage{amsmath}
\usepackage{amsfonts}
\usepackage{amssymb}
\usepackage{hyperref}

\usepackage{graphicx}
\providecommand{\U}[1]{\protect\rule{.1in}{.1in}}

\begin{document}

\title{Groups of order $p^{9}$, class 2, and exponent $p$ having derived group of order $p^{2}$}
\author{Doug Tyler}
\date{October 2017}   
\maketitle

\begin{abstract}

This paper concerns finite groups of class (at most) two and of odd prime exponent $p$. Such a group is called \textit{special} if the center lies within its derived group,  $Z(G) \leq G'$ (the general definition of \textit{special}, without the exponent $p$ assumption, is more involved). Every group of class 2 and exponent $p$ can be uniquely expressed as the direct product of an elementary abelian group and a special group. This reduces the isomorphism problem to special groups. The special groups having $G'$ cyclic are well known. Groups having $|G'|=p^{2}$ are known to be generated by two abelian subgroups. As such, they can be described by a pair of  \textit{Scharlau Matrices} which we will define. Using these, Vishnevetskii ([1], [2]) classified the special groups which are not central products of groups of smaller order. We call these  \textit{Vishnevetskii indecomposable}. All decomposable groups are central products of two or more indecomposable groups.  By Theorem 2 of [2], if $G$ is a special group with derived group of order $p^2$ then the indecomposable central factors of $G$, together with their multiplicities, form a set of invariants for $G$. However these invariants do not determine $G$, since (as we will show below) two nonisomorphic groups can have the same indecomposable central factors. The groups of exponent $p$ and order dividing $p^{8}$ are already in the literature. In this paper we prove that there are six isomorphism types of the special groups of order $p^{9}$ having $|G'|=p^{2}$, and we list them. Five of these groups are central products of two or more indecomposable groups. One of them, which Vishnevetskii described, is indecomposable and is not a central product of smaller groups.

We will also introduce a way of representing groups of class 2 having exponent $p$ by using digraphs with flows on the edges. The digraph is compact and gives, on sight, a great deal of structural information about the group and completely defines the group. Some group invariants will also be described which are easy to compute from the digraphs and which can readily be used to distinguish the isomorphism types for small orders. 

For each of the special groups of order dividing $p^{9}$ we will list a pair of Scharlau matrices defining the group, we will identify the central factors when they exist, and we will give group invariants. For the groups of orders dividing $p^{8}$ we identify where their power commutator presentations appear in the literature. For the groups of order $p^{9}$ we will give presentations as digraphs. 

\end{abstract}

\section{Background}                                 

If $G$ has class at most 2 then its derived group lies in its center, $G'\le Z(G)$. If $G$ has exponent $p$ then $G'$ and the Frattini subgroup, $\Phi(G)$, coincide. It is easy to show that if $G$ has exponent $p$ then $G$ has a cyclic direct factor if and only if $G'<Z(G)$. Thus, to classify groups of class 2 and exponent $p$, we only have to treat those with $G'=Z(G)=\Phi(G)$. These groups are called  \textit{special}. If $G$ is special and $|G'|=p$ then $G$ is called  \textit{extraspecial}. We denote the elementary abelian groups by $Z_p^n$ for $n=1,2,...$ and note that these are the groups of class 1 or, equivalently, the groups with $G'=1$. The extraspecial groups are also known, there being exactly one group each (of exponent $p$) of sizes  $p^{3},p^{5},p^{7},...$. We will denote these by $E_3,E_5,E_7,...$. The group $E_3$ is given by %
$E_3=\langle x_1,x_2\,|\,x_1^p=1,x_2^p=1,[x_1,x_2]^p=1,\,\text{class}\,2\rangle$. The central product of $n$ copies of $E_3$ can be taken in only one way and is isomorphic to $E_{2n+1}$. Thus if $G$ has exponent $p$ and $|G'|=p$ then $G$ has the form $Z_p^m\times E_{2n+1}$.

For groups of class 2 having exponent $p$,  $G/G'$ and $G'$ are elementary abelian and $G$ is completely determined by how the commutators of a generating set of $G$ lie inside $G'$. If $|G/G'|=p^d$ and $|G'|=p^r$ with $G=\langle x_1,x_2,...,x_d \rangle$ and $G'=\langle z_1,z_2,...,z_r\rangle $  then $G$ is determined by the values of the $a_k(i,j)$'s in the relations%
\[
[x_i,x_j]=z_1^{a_1(i,j)}z_2^{a_2(i,j)}...z_r^{a_r(i,j)}, \, \text{for}\, 1\le i<j\le d.
\]
The undefined entries can be filled in using $a_k(j,i)=-a_k(i,j)$ for $j > i$. For example, the smallest group of exponent $p$, having class 2 and $|G'|=p^2$ is the group of order $p^5$ with $G=\langle x_1,x_2,x_3\rangle$, $G'=\langle z_1,z_2\rangle$, $[x_1,x_3]=1$, $[x_1,x_2]=z_1$, and $[x_2,x_3]=z_2$. Expressed in terms of generators and relations this is%
\[
G=\langle x_1,x_2,x_3\,|\,[x_1,x_3]=1,\, \text{exponent}\,p,\,\text{class}\,2\rangle .
\]
This tells us $[x_1,x_3]=1$ and that $[x_1,x_2]$ and $[x_2,x_3]$ are independent elements of $G'$. The 3 by 3 antisymmetric matrices,  $A_k=((a_k(i,j))_{i,j})$ for $k=1,2$ are easily worked out. 

When $|G'|=p^2$, the array, in this case the pair, of antisymmetric matrices takes on a particularly simple form. Specifically,  $G$ necessarily has two elementary abelian subgroups%
\[
X=\langle x_1,x_2,...,x_m,z_1,z_2\rangle  \, \text{and}\, Y=\langle y_1,y_2,...,y_n,z_1,z_2\rangle
\]
such that $G=XY$, $G'=X\cap Y$, $m+n=d$, and  $m$ and $n$  differ by at most 1. (See Scharlau [3] and Sergeichuk [4].) Because of this, there are $m \times n$ matrices $A=(a(i,j))$ and $B=(b(i,j))$ (which are smaller and are not antisymmetric) which we call  \textit{Scharlau matrices} such that%
\[
[x_i,y_j]=z_1^{a(i,j)}z_2^{b(i,j)},\, \text{for} \, 1\le i\le m,\, 1\le j\le n.  
\]
The isomorphism type of the group is determined by the matrices $A$ and $B$. (Given $G$ there will generally be multiple ways to choose $A$ and $B$.) Vishnevetskii ([1], [2]) investigated the possible pairs of Scharlau matrices using earlier work of Kronecker and Dieudonne on matrices. We will say that $G$ is \textit{Vishnevetskii indecomposable} if $G$ is not a central product of smaller groups. Vishnevetskii classified those groups which are special and indecomposable, producing a list of the distinct isomorphism types. All special groups which are decomposable are central products of smaller indecomposables. According to Vishnevetskii's Theorem 2 in [2], the isomorphism types, and multiplicities, of the central factors are an invariant for a decomposable group. We shall see in Section 3 that the groups 8.6.7 and 8.6.8 have the same central factors and multiplicities but are not isomorphic. This holds for groups 8.6.9 and 8.6.10 also. In Section 4 the pair of groups named ``B" and ``C" provide an example of size $p^9$.

A contemporary list of groups of exponent $p$ having order dividing $p^8$, using power commutator presentations, can be found at Vaughan-Lee ([6]) We will reproduce Section 3 those groups which are special and have $|G'|=p^2$, giving Scharlau matrices, identifying the central factors, and tabulating certain invariants which we describe in Section 2. We list these groups in the order in which they appear in [6]. (The frequencies of central quotients are sufficient to distinguish all but two of the groups.)

We do not report here what the indecomposable groups are that Vishnevetskii developed. A succinct description of the Scharlau matrices has recently been given by Vaughan-Lee in the first two pages of [5]. The indecomposable groups having orders dividing $p^9$ are listed in this paper in the sections below. The names are 5.3.1, 6.4.3, 6.4.4, 7.5.6, 8.6.13, and 8.6.14 for those having order dividing $p^8$ and the group named ``F" in the sections on groups of order $p^9$.

\section{Digraphs and invariants of groups of exponent $p$ having class 2}                 

Given a group of exponent $p$ and class 2, the usual way of defining it is to give its commutator presentation on a set of generators, $G=\langle x_1,x_2,...,x_d \rangle $.  The size of $G'$ must be inferred from the relations on the generating set. 

We introduce here a way to represent, and define, such a group by a digraph with flows on the edges. It is necessary to assign names to the elements of a basis for $G'$, $\langle z_1,z_2,...,z_r \rangle$. The digraph of the group has vertices labelled with the generators of $G$. For each pair of generators  $(x_i,x_j)$ with $i\le j$, the directed edge from $x_i$ to $x_j$ is labelled with a flow equal to the element of $G'$ corresponding to the value of their commutator, $[x_i,x_j]$, expressed in terms of the $z_k$'s. (The direction of the edge can be reversed if the flow is replaced by its inverse, $[x_j,x_i]$.) For generators which commute the flow is the identity element of $G'$ and the edge is omitted. Isolated vertices with no edges correspond to cyclic direct factors and will be absent if and only if the group is special. A representative extraspecial group is $E_5$. Its digraph looks like
\\
\setlength{\unitlength}{1.0cm}
\begin{picture}(9,1.0)
\thicklines
\put(4.0,0.1){$x_1$}
\put(5.6,0.1){$x_2$}
\put(7.2,0.1){$x_3$}
\put(8.8,0.1){$x_4$}
\put(4.8,0.5){$z_1$}
\put(8.0,0.5){$z_1$}
\put(4.5,0.2){\vector(1,0){1.0}}
\put(7.7,0.2){\vector(1,0){1.0}}
\end{picture}
\\

\noindent
Notice that the group is the central product of two copies of $E_3$ and that these are the subgroups generated by the components of the graph. 

The group of order $p^5$ defined in the previous section is completely defined by the digraph 
\\
\setlength{\unitlength}{1.0cm}
\begin{picture}(10,1.0)
\thicklines
\put(4.5,0.1){$x_1$}
\put(6.3,0.1){$x_2$}
\put(8.1,0.1){$x_3$}
\put(5.4,0.6){$z_1$}
\put(7.2,0.6){$z_2$}
\put(5.1,0.2){\vector(1,0){1.0}}
\put(6.9,0.2){\vector(1,0){1.0}}
\end{picture}
\\

\noindent
Digraphs cost us the inconvenience of having to name the generators of $G'$ but buy us immediate information about the size of $G'$ and a view of commuting pairs of generators. Of greater value is the fact that the impact of changes of variables in each of (the vector spaces) $G/G'$ and $G'$ can be computed on sight. For example replacing $x_2$ with $x_1^a x_2$ in the picture above has no impact on the flow as $x_1$ already commutes with itself. Replacing $x_3$ with $x_1^a x_3^b$ where $b\ne 0$ will change the flow $z_2$  to $z_1^{-a} z_2^b$. Replacing $x_1$  with $x_1 x_2^a$ will introduce an edge, with a flow, from $x_1$ to $x_3$. The centralizers of the elements of $G$ are obvious --- the elements inside $\langle x_1,x_3,G' \rangle$ have centralizers of index at most $p$ (and generate a copy of $Z_p^4$) while the elements outside have centralizers of order $p^3$. From the digraph one can see that the automorphism group of this group acts like the general linear group GL(2,$p$) on its derived group, a fact we will use later.

Using the transformations on the digraphs provides a particularly efficient way to prove the classification of the extraspecial groups. For groups with $|G'|=p^r \ge p^2$, it is easy to compute the isomorphism type (and center) of $G/\langle z \rangle$ for any cyclic, central subgroup $\langle z \rangle \ne 1$. Since $\langle z \rangle \subset G'=Z(G)$, the quotient has exactly the same number of generators as $G$ and has commutator subgroup of size exactly $p^{r-1}$. The quotient will not, in general, be special and the center may be larger than $p^{r-1}$. There are precisely $(p^r-1)/(p-1)$  of these quotients and the resulting isomorphism types, counted according to multiplicity, form a group invariant. In the special group of order $p^5$ above $G/ \langle z \rangle$ is isomorphic to $Z_p \times E_3$ for all $p+1$ quotients. 

For groups with$|G'|=p^2$, there are exactly $p+1$ choices of $\langle z \rangle$ and exactly $p+1$ central quotients. They will all be of the form $Z_p^{d-2n} \times E_{2n+1}$ and all that is needed is to count the number of occurrences of $E_{2n+1}$ for $n=1,2,..., \lfloor d/2 \rfloor$. This is easily tabulated in a vector of length $\lfloor d/2 \rfloor$ and is done in the tables in the sections below for each of the special groups of order dividing $p^9$. Thus, for example, the first entry in the table of groups of orders dividing $p^8$ below, is the group of order $p^5$. Its central quotient invariant is the ``vector" of length 1 telling us that $E_3$ appears exactly $p+1$ times. 

\section{The special groups of orders dividing $p^8$ having derived group of order $p^2$}      

We list here the isomorphism types of groups of order at most $p^8$, having exponent $p$, class 2 and $|G'|=p^2$ . We follow the notation of Vaughan-Lee ([6]) but limit our attention to special groups. Each group is denoted by three numbers, the first of which indicates the order and the second of which denotes the size of a minimal generating set. The third number is an index for the isomorphism type. Vaughan-Lee describes a scheme for ordering the isomorphism types which works up to this size. Note that the group $E_3$ is the ONLY Vishnevetskii indecomposable group which has $|G'| \le p$. It appears almost immediately as a central factor in larger groups and is group 3.2.1 in [6]. Its digraph is the first component of the digraph of $E_5$ shown in the previous section.

We tabulate the groups of order below $p^8$ first. For groups which are decomposable, the isomorphism types of the central factors are listed. Matrices $A$ and $B$ are the Scharlau Matrices and provide defining relations. Finally, the isomorphism types of the central quotients are counted. (Digraphs of the indecomposable groups of orders below $p^8$ can be found in the components of the special groups of order $p^9$ in the next section.)
\\

\begin{center}
\begin{tabular}{ | c || c | c | c | c | c | } 
\hline
$\begin{array} [c]{c} \text{Group \#}    \\ \text{in [6]} \end{array}$ &
$\begin{array} [c]{c} \text{Order}         \\  \end{array}$ &
$\begin{array} [c]{c} \text{Central}       \\ \text{Factors} \end{array}$ &
$\begin{array} [c]{c} \text{Scharlau}     \\ \text{Matrix~} A \end{array}$ &
$\begin{array} [c]{c} \text{Scharlau}     \\ \text{Matrix~} B \end{array}$ &
$\begin{array} [c]{c} G/\langle z\rangle \\ \text{freq.} \end{array}$ \\ \hline \hline

5.3.1 & $p^5$ & indec. & 
$ \left( \begin{array} [c]{cc} 1 & 0 \end{array} \right) $ &
$ \left( \begin{array} [c]{cc} 0 & 1 \end{array} \right) $ & $(p+1) $ \\ \hline

6.4.2 & $p^6$ & $ \begin{array} [c]{c} E_3 \\ E_3 \end{array} $ &
$ \left( \begin{array} [c]{cc} 1 & 0\\ 0 & 0 \end{array} \right) $ & 
$ \left( \begin{array} [c]{cc} 0 & 0\\ 0 & 1 \end{array} \right) $ & $(2,p-1)$\\ \hline

6.4.3 & $p^6$ & indec. &
$ \left( \begin{array} [c]{cc} 1 & 0\\ 0 & 1 \end{array} \right) $ &
$ \left( \begin{array} [c]{cc} 0 & 1\\ 0 & 0 \end{array} \right) $ & $(1,p)$\\ \hline

6.4.4 & $p^6$ & indec. &
$ \left( \begin{array} [c]{cc} 1 & 0\\ 0 & 1 \end{array} \right) $ &
$ \left( \begin{array} [c]{cc} 0 & 1\\ \nu & 0 \end{array} \right) $ & $(0,p+1)$\\ \hline

7.5.5 & $p^7$ & $ \begin{array} [c]{c} E_3\\ 5.3.1 \end{array} $ &
$ \left( \begin{array} [c]{ccc} 1 & 0 & 0\\ 0 & 0 & 1 \end{array} \right) $ &
$ \left( \begin{array} [c]{ccc} 0 & 1 & 0\\ 0 & 0 & 0 \end{array} \right) $ & 
$(1,p)$\\ \hline

7.5.6 & $p^7$ & indec. &
$ \left( \begin{array} [c]{ccc} 1 & 0 & 0\\ 0 & 1 & 0 \end{array} \right) $ &
$ \left( \begin{array} [c]{ccc} 0 & 1 & 0\\ 0 & 0 & 1 \end{array} \right) $ &
$(0,p+1)$\\ \hline
\end{tabular}
\end{center}

\begin{center}
\bf Table of the six (6) special groups of orders below $p^8$
\end{center}

In group 6.4.4, $\nu$ denotes an arbitrary quadratic nonresidue modulo $p$. Thus Matrix $B$ denotes the companion matrix to the (arbitrary) irreducible quadratic, $x^2-\nu$,  modulo $p$. (This applies to group 8.6.11 below also.) The groups of orders below $p^8$  can be distinguished by their central quotients.

We now tabulate the groups of order $p^8$. 
\\

\begin{center}
\begin{tabular}{ | c || c | c | c | c | } 
\hline
$\begin{array} [c]{c} \text{Group \#}    \\ \text{in [6]} \end{array}$ &
$\begin{array} [c]{c} \text{Central}       \\ \text{Factors} \end{array}$ &
$\begin{array} [c]{c} \text{Scharlau}     \\ \text{Matrix~} A \end{array}$ &
$\begin{array} [c]{c} \text{Scharlau}     \\ \text{Matrix~} B \end{array}$ &
$\begin{array} [c]{c} G/\langle z\rangle \\ \text{frequencies} \end{array}$ \\ \hline \hline

8.6.7  & $ \begin{array} [c]{c} E_3\\ E_3\\ E_3 \end{array} $ & 
$ \left( \begin{array} [c]{ccc} 1 & 0 & 0\\ 0 & 1 & 0\\ 0 & 0 & 0 \end{array} \right) $ &
$ \left( \begin{array} [c]{ccc} 0 & 0 & 0\\ 0 & 0 & 0\\ 0 & 0 & 1 \end{array} \right) $ & 
$(1,1,p-1) $ \\ \hline

8.6.8  & $ \begin{array} [c]{c} E_3\\ E_3\\ E_3 \end{array} $ &
$ \left( \begin{array} [c]{ccc} 1 & 0 & 0\\ 0 & 1 & 0\\ 0 & 0 & 0 \end{array} \right) $ & 
$ \left( \begin{array} [c]{ccc} 0 & 0 & 0\\ 0 & 1 & 0\\ 0 & 0 & 1 \end{array} \right) $ & 
$(0,3,p-2)$\\ \hline

8.6.9  & $ \begin{array} [c]{c} E_3\\ 6.4.3 \end{array} $ &
$ \left( \begin{array} [c]{ccc} 1 & 0 & 0\\ 0 & 1 & 0\\ 0 & 0 & 1 \end{array} \right) $ & 
$ \left( \begin{array} [c]{ccc} 0 & 1 & 0\\ 0 & 0 & 0\\ 0 & 0 & 1 \end{array} \right) $ & 
$(0,2,p-1)$\\ \hline

8.6.10  & $ \begin{array} [c]{c} E_3\\ 6.4.3 \end{array} $ &
$ \left( \begin{array} [c]{ccc} 1 & 0 & 0\\ 0 & 1 & 0\\ 0 & 0 & 0 \end{array} \right) $ & 
$ \left( \begin{array} [c]{ccc} 0 & 1 & 0\\ 0 & 0 & 0\\ 0 & 0 & 0 \end{array} \right) $ & 
$(1,0,p)$\\ \hline

8.6.11 & $ \begin{array} [c]{c} E_3\\ 6.4.4 \end{array} $ &
$ \left( \begin{array} [c]{ccc} 1 & 0 & 0\\ 0 & 1 & 0\\ 0 & 0 & 1 \end{array} \right) $ & 
$ \left( \begin{array} [c]{ccc} 0 & 1 & 0\\ \nu & 0 & 0\\ 0 & 0 & 0 \end{array} \right) $ & 
$(0,1,p)$\\ \hline

8.6.12 & $ \begin{array} [c]{c} 5.3.1\\ 5.3.1 \end{array} $ &
$ \left( \begin{array} [c]{ccc} 1 & 0 & 0\\ 0 & 0 & 1\\ 0 & 0 & 0 \end{array} \right) $ & 
$ \left( \begin{array} [c]{ccc} 0 & 1 & 0\\ 0 & 0 & 0\\ 0 & 0 & 1 \end{array} \right) $ & 
$(0,p+1,0)$\\ \hline

8.6.13 & indec. &
$ \left( \begin{array} [c]{ccc} 1 & 0 & 0\\ 0 & 1 & 0\\ 0 & 0 & 1 \end{array} \right) $ & 
$ \left( \begin{array} [c]{ccc} 0 & 1 & 0\\ 0 & 0 & 1\\ 0 & 0 & 0 \end{array} \right) $ & 
$(0,1,p)$\\ \hline

8.6.14 & indec. &
$ \left( \begin{array} [c]{ccc} 1 & 0 & 0\\ 0 & 1 & 0\\ 0 & 0 & 1 \end{array} \right) $ & 
$ \left( \begin{array} [c]{ccc} 0 & 1 & 0\\ 0 & 0 & 1\\ c & b & a \end{array} \right) $ & 
$(0,0,p+1)$\\ \hline
\end{tabular}
\end{center}

\begin{center}
\bf Table of the eight (8) special groups of order $p^8$
\end{center}

In group 8.6.14, Matrix $B$ denotes the companion matrix to an arbitrary irreducible cubic modulo $p$. (Different irreducible cubics give isomorphic groups.) The eight groups of order $p^8$ can be distinguished by their central quotients with the exception of groups 8.6.11 and 8.6.13. These two can be distinguished by the fact that 8.6.13 is indecomposable while 8.6.11 is not. Alternatively, the elements with centralizers of index at most $p$ commute in 8.6.13 but do not in 8.6.11.

Note also that a decomposable group is not, in general, determined by its central factors. Groups 8.6.7 and 8.6.8 are central products of three copies of group $E_3$ yet they are not isomorphic. Similarly a central product of $E_3$ and 6.4.3 can be either 8.6.9 or 8.6.10 and these are not isomorphic.

\section{There are at least six (special) groups of order $p^9$ having $|G'|=p^2$}     

      This section tabulates six groups called A, B, C, D, E, and F. The table has the central factors, the Scharlau matrices and the frequencies of the central quotients. 

The central quotient frequencies distinguish between the groups, except for D and F and for A and E. For group D the elements with centralizers of index at most $p$ form a subgroup and for group F they do not. Thus D and F cannot be isomorphic. For group E, if $G/ \langle z \rangle$ is not extraspecial then the preimage of the center of $G/ \langle z \rangle$ is an abelian subgroup of E. For group A the resulting subgroup is nonabelian so that A and E are not isomorphic.

An entirely different argument comes from the fact that the indecomposable central factors of the six groups are different with the exception of groups B and C. Vishnevetskii stated that the sizes and isomorphism types of the central factors are a group invariant. This distinguishes all of the groups except B and C. The central quotients of B and C are different so they are not isomorphic.

Note that groups B and C provide another example of non-isomorphic groups which have central factors which are the same up to isomorphism and multiplicities.
\\
We first illustrate the digraphs of these groups. Note that the components of the digraphs indicate what the central factors are.
\\
A. 	
\setlength{\unitlength}{1.0cm}
\begin{picture}(11,1.0)
\thicklines
\put(0.4,0.1){$x_1$}
\put(2.0,0.1){$x_2$}
\put(3.6,0.1){$x_3$}
\put(5.2,0.1){$x_4$}
\put(6.8,0.1){$x_5$}
\put(8.4,0.1){$x_6$}
\put(10.0,0.1){$x_7$}
\put(1.2,0.6){$z_1$}
\put(2.8,0.6){$z_2$}
\put(4.4,0.6){$z_1$}
\put(6.0,0.6){$z_2$}
\put(9.2,0.6){$z_1$}
\put(0.9,0.2){\vector(1,0){1.0}}
\put(2.5,0.2){\vector(1,0){1.0}}
\put(4.1,0.2){\vector(1,0){1.0}}
\put(5.7,0.2){\vector(1,0){1.0}}
\put(8.9,0.2){\vector(1,0){1.0}}
\end{picture}
\\
B. 	
\setlength{\unitlength}{1.0cm}
\begin{picture}(11,1.0)
\thicklines
\put(0.4,0.1){$x_1$}
\put(2.0,0.1){$x_2$}
\put(3.6,0.1){$x_3$}
\put(5.2,0.1){$x_4$}
\put(6.8,0.1){$x_5$}
\put(8.4,0.1){$x_6$}
\put(10.0,0.1){$x_7$}
\put(1.2,0.6){$z_1$}
\put(2.8,0.6){$z_2$}
\put(6.0,0.6){$z_1$}
\put(9.2,0.6){$z_1$}
\put(0.9,0.2){\vector(1,0){1.0}}
\put(2.5,0.2){\vector(1,0){1.0}}
\put(5.7,0.2){\vector(1,0){1.0}}
\put(8.9,0.2){\vector(1,0){1.0}}
\end{picture}
\\
C. 
\setlength{\unitlength}{1.0cm}
\begin{picture}(11,1.0)
\thicklines
\put(0.4,0.1){$x_1$}
\put(2.0,0.1){$x_2$}
\put(3.6,0.1){$x_3$}
\put(5.2,0.1){$x_4$}
\put(6.8,0.1){$x_5$}
\put(8.4,0.1){$x_6$}
\put(10.0,0.1){$x_7$}
\put(1.2,0.6){$z_1$}
\put(2.8,0.6){$z_2$}
\put(6.0,0.6){$z_1$}
\put(9.2,0.6){$z_2$}
\put(0.9,0.2){\vector(1,0){1.0}}
\put(2.5,0.2){\vector(1,0){1.0}}
\put(5.7,0.2){\vector(1,0){1.0}}
\put(8.9,0.2){\vector(1,0){1.0}}
\end{picture}
\\
D. 
\setlength{\unitlength}{1.0cm}
\begin{picture}(11,2.0)
\thicklines
\put(0.9,0.9){$x_1$}
\put(2.5,1.6){$x_2$}
\put(2.5,0.2){$x_3$}
\put(4.1,0.9){$x_4$}
\put(6.1,0.9){$x_5$}
\put(7.7,0.9){$x_6$}
\put(9.3,0.9){$x_7$}
\put(1.5,1.4){$z_1$}
\put(1.6,0.4){$z_2$}
\put(3.5,0.3){$z_1^\nu$}
\put(3.5,1.5){$z_2$}
\put(6.9,1.4){$z_1$}
\put(8.5,1.4){$z_2$}
\put(1.4,1.1){\vector(2,1){1.0}}
\put(1.4,0.9){\vector(2,-1){1.0}}
\put(3.0,1.6){\vector(2,-1){1.0}}
\put(3.0,0.4){\vector(2,1){1.0}}
\put(6.6,1.0){\vector(1,0){1.0}}
\put(8.2,1.0){\vector(1,0){1.0}}
\end{picture}	
\\   
E. 
\setlength{\unitlength}{1.0cm}
\begin{picture}(11,1.0)
\thicklines
\put(0.4,0.1){$x_1$}
\put(2.0,0.1){$x_2$}
\put(3.6,0.1){$x_3$}
\put(5.2,0.1){$x_4$}
\put(6.8,0.1){$x_5$}
\put(8.4,0.1){$x_6$}
\put(10.0,0.1){$x_7$}
\put(1.2,0.6){$z_1$}
\put(2.8,0.6){$z_2$}
\put(4.4,0.6){$z_1$}
\put(7.6,0.6){$z_1$}
\put(9.2,0.6){$z_2$}
\put(0.9,0.2){\vector(1,0){1.0}}
\put(2.5,0.2){\vector(1,0){1.0}}
\put(4.1,0.2){\vector(1,0){1.0}}
\put(7.3,0.2){\vector(1,0){1.0}}
\put(8.9,0.2){\vector(1,0){1.0}}
\end{picture}	
\\
F. 
\setlength{\unitlength}{1.0cm}
\begin{picture}(11,1.0)
\thicklines
\put(0.4,0.1){$x_1$}
\put(2.0,0.1){$x_2$}
\put(3.6,0.1){$x_3$}
\put(5.2,0.1){$x_4$}
\put(6.8,0.1){$x_5$}
\put(8.4,0.1){$x_6$}
\put(10.0,0.1){$x_7$}
\put(1.2,0.6){$z_1$}
\put(2.8,0.6){$z_2$}
\put(4.4,0.6){$z_1$}
\put(6.0,0.6){$z_2$}
\put(7.6,0.6){$z_1$}
\put(9.2,0.6){$z_2$}
\put(0.9,0.2){\vector(1,0){1.0}}
\put(2.5,0.2){\vector(1,0){1.0}}
\put(4.1,0.2){\vector(1,0){1.0}}
\put(5.7,0.2){\vector(1,0){1.0}}
\put(7.3,0.2){\vector(1,0){1.0}}
\put(8.9,0.2){\vector(1,0){1.0}}
\end{picture}
\\
\begin{center}
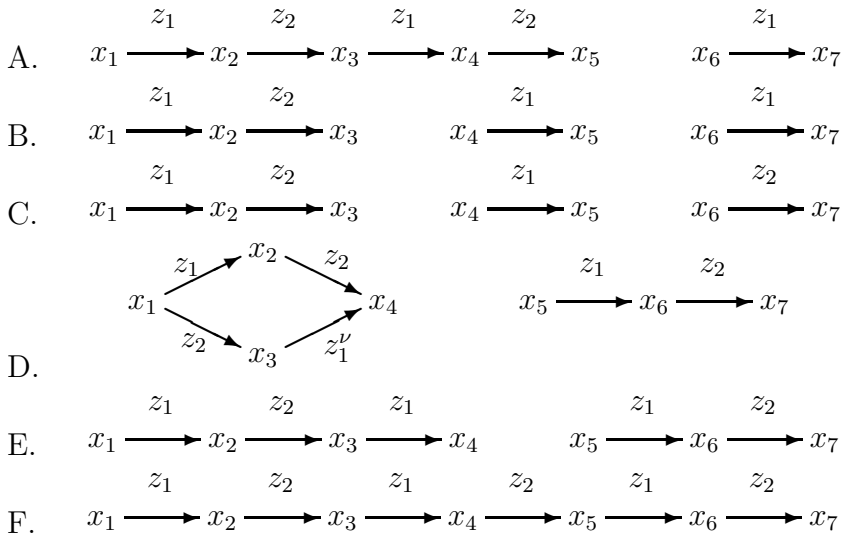
\bf Digraphs of the Six Special Groups of Order $p^9$
\end{center}

\begin{center}
\begin{tabular}{ | c || c | c | c | c | } 
\hline
$\begin{array} [c]{c} \text{Group}    \\ \text{name} \end{array}$ &
$\begin{array} [c]{c} \text{Central}       \\ \text{Factors} \end{array}$ &
$\begin{array} [c]{c} \text{Scharlau}     \\ \text{Matrix~} A \end{array}$ &
$\begin{array} [c]{c} \text{Scharlau}     \\ \text{Matrix~} B \end{array}$ &
$\begin{array} [c]{c} G/\langle z\rangle \\ \text{frequencies} \end{array}$ \\ \hline \hline

A  & $ \begin{array} [c]{c} E_3\\ 7.5.6 \end{array} $ &
$ \left( \begin{array} [c]{cccc} 1 & 0 & 0 & 0\\ 0 & 1 & 0 & 0\\ 0 & 0 & 1 & 0 \end{array} \right) $
& 
$ \left( \begin{array} [c]{cccc} 0 & 0 & 0 & 0\\ 0 & 0 & 1 & 0\\ 0 & 0 & 0 & 1 \end{array} \right) $
& $(0,1,p)$\\ \hline

B  & $ \begin{array} [c]{c} E_3\\ E_3\\ 5.3.1 \end{array} $ &
$ \left( \begin{array} [c]{cccc} 1 & 0 & 0 & 0\\ 0 & 1 & 0 & 0\\ 0 & 0 & 1 & 0 \end{array} \right) $
& 
$ \left( \begin{array} [c]{cccc} 0 & 0 & 0 & 0\\ 0 & 0 & 0 & 0\\ 0 & 0 & 0 & 1 \end{array} \right) $
& $(1,0,p)$\\ \hline

C  & $ \begin{array} [c]{c} E_3\\ E_3\\ 5.3.1 \end{array} $ &
$ \left( \begin{array} [c]{cccc} 1 & 0 & 0 & 0\\ 0 & 0 & 0 & 0\\ 0 & 0 & 1 & 0 \end{array} \right) $
& 
$ \left( \begin{array} [c]{cccc} 0 & 0 & 0 & 0\\ 0 & 1 & 0 & 0\\ 0 & 0 & 0 & 1 \end{array} \right) $
& $(0,2,p-1)$\\ \hline

D & $ \begin{array} [c]{c} 5.3.1\\ 6.4.4 \end{array} $ &
$ \left( \begin{array} [c]{cccc} 1 & 0 & 0 & 0\\ 0 & 0 & 1 & 0\\ 0 & 0 & 0 & 1 \end{array} \right) $
& 
$ \left( \begin{array} [c]{cccc} 0 & 1 & 0 & 0\\ 0 & 0 & 0 & 1\\ 0 & 0 & \nu & 0 \end{array} \right) $
& $(0,0,p+1)$\\ \hline

E & $ \begin{array} [c]{c} 5.3.1\\ 6.4.3 \end{array} $ &
$ \left( \begin{array} [c]{cccc} 1 & 0 & 0 & 0\\ 0 & 1 & 0 & 0\\ 0 & 0 & 1 & 0 \end{array} \right) $
& 
$ \left( \begin{array} [c]{cccc} 0 & 1 & 0 & 0\\ 0 & 0 & 0 & 0\\ 0 & 0 & 0 & 1 \end{array} \right) $
& $(0,1,p)$\\ \hline

F & indec. &
$ \left( \begin{array} [c]{cccc} 1 & 0 & 0 & 0\\ 0 & 1 & 0 & 0\\ 0 & 0 & 1 & 0 \end{array} \right) $
& 
$ \left( \begin{array} [c]{cccc} 0 & 1 & 0 & 0\\ 0 & 0 & 1 & 0\\ 0 & 0 & 0 & 1 \end{array} \right) $
& $(0,0,p+1)$\\ \hline
\end{tabular}
\end{center}

\begin{center}
\bf Table of the six (6) special groups of order $p^9$
\end{center}

\section{There are at most six (special) groups of order $p^9$ having $|G'|=p^2$}     

In Vishnevetskii's work on indecomposable groups, he proves that there is at most one such group of order $p^{2n+1}$ for any $n$. For $p^9$ this is the group called ``F" in the previous section. The remaining groups must be central products of smaller groups and, with the exception of the factor $E_3$, all factors must be special and indecomposable and have derived group of order $p^2$. The only other possible factors are therefore 5.3.1, 6.4.3, 6.4.4 and 7.5.6. The digraph of $G$ will have seven vertices and the vertices in the various components must come from the central factors. The groups 5.3.1, 6.4.3, 6.4.4 and 7.5.6 contribute 3, 4, 4 and 5 vertices to $G$ and the group $E_3$ contributes two vertices. The central factors therefore provide a partition of the number 7 into parts of size at least 2. The only possibilities are 7 (the indecomposable case), 5+2, 4+3 and 3+2+2. 

Since the only group contributing three vertices is 5.3.1, the partitions involving 3 have 5.3.1 as a central factor. The partition 3+2+2 corresponds to the central product of three factors, two isomorphic to $E_3$ and the third isomorphic to 5.3.1. Also, the flows on the edges of the two copies of $E_3$ either generate $G'$ or can be taken to be equal. The partition 4+3 must be a central product of 5.3.1 with either 6.4.3 or 6.4.4. Finally, since the automorphism group of 5.3.1 acts like the general linear group, GL(2,$p$), on its derived group, the flows on the edges of the group 5.3.1 can be chosen to equal $z_1$ and $z_2$ as needed. Suppose, for example, that we have a central product of 6.4.4 and 5.3.1. The digraph of 6.4.4 is given by the component on the left in the diagram for group D above and we can choose $z_1$ and $z_2$ from the derived group to match these edges in 6.4.4. The central product must have the form of the digraph
\\
\setlength{\unitlength}{1.0cm}
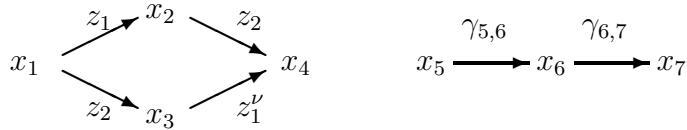
\begin{picture}(11,2.0)
\thicklines
\put(1.7,0.9){$x_1$}
\put(3.5,1.6){$x_2$}
\put(3.5,0.2){$x_3$}
\put(5.3,0.9){$x_4$}
\put(7.1,0.9){$x_5$}
\put(8.7,0.9){$x_6$}
\put(10.3,0.9){$x_7$}
\put(2.7,1.5){$z_1$}
\put(2.7,0.3){$z_2$}
\put(4.7,0.3){$z_1^\nu$}
\put(4.7,1.5){$z_2$}
\put(7.7,1.4){$ \gamma_{5,6} $}
\put(9.3,1.4){$ \gamma_{6,7} $}
\put(2.4,1.1){\vector(2,1){1.0}}
\put(2.4,0.9){\vector(2,-1){1.0}}
\put(4.1,1.6){\vector(2,-1){1.0}}
\put(4.1,0.4){\vector(2,1){1.0}}
\put(7.6,1.0){\vector(1,0){1.0}}
\put(9.2,1.0){\vector(1,0){1.0}}
\end{picture}	
\\   
Because of the GL(2,$p$) action on the derived group of 5.3.1, there is a change of variables that can be made with $x_5$ and $x_7$ so that the new edges $ \gamma_{5,6}$ and  $ \gamma_{6,7}$  become $z_1$ and $z_2$. Thus this group must be isomorphic to group D. A completely analogous argument can be applied to show that group E is the only possible result of a central product of 6.4.3 and 5.3.1. Finally, suppose we have a central product of 5.3.1 with two copies of $E_3$. If the derived groups of the two $E_3$'s do not coincide then we can label the edges with the commutators $z_1$ and $z_2$. Arguing as before, the action on the derived group of 5.3.1 enables us to take the edges to these two elements, giving us group C. If the derived groups of the two $E_3$'s coincide, we may label one edge with $z_1$ and then make a change of variables of the form $x^\prime = x^\alpha$ to carry the next edge to $z_1$. Finally, we can again change variables in the generators of the copy of 5.3.1 to cause one of it's edges to equal $z_1$. The last edge of 5.3.1 can be taken to be $z_2$ and we get group B. Thus, the groups involving 6.4.3 and 6.4.4 are E and D. The other two groups are B and C (where the flows from the $E_3$'s do not, or do, generate all of $G'$).

The partition 5+2 must involve the central product of 7.5.6 with $E_3$. One such possibility is group A. We use digraphs to show that group A is the only isomorphism type here. Certainly any such central product must have the form 
\\

\setlength{\unitlength}{1.0cm}
\begin{picture}(12,1.0)
\thicklines
\put(1.4,0.1){$x_1$}
\put(3.0,0.1){$x_2$}
\put(4.6,0.1){$x_3$}
\put(6.2,0.1){$x_4$}
\put(7.8,0.1){$x_5$}
\put(9.4,0.1){$x_6$}
\put(11.0,0.1){$x_7$}
\put(2.2,0.6){$z_1$}
\put(3.8,0.6){$z_2$}
\put(5.4,0.6){$z_1$}
\put(7.0,0.6){$z_2$}
\put(10.2,0.6){$z$}
\put(1.9,0.2){\vector(1,0){1.0}}
\put(3.5,0.2){\vector(1,0){1.0}}
\put(5.1,0.2){\vector(1,0){1.0}}
\put(6.7,0.2){\vector(1,0){1.0}}
\put(9.9,0.2){\vector(1,0){1.0}}
\end{picture}
\\

\noindent
where $z=z_1^a z_2^b \ne 1$ . From the mirror symmetry in $\langle x_1,x_2,x_3,x_4,x_5 \rangle$ (group 7.5.6), we may interchange a and b. We may as well assume that $a \ne 0$. Changes of variables of the form $x_i^\prime = x_i^{\alpha (i)} ,\; z_i^\prime = z_i^{\beta (i)}$ ,  allow us to assume that $a = 1$ and that $b$ is either 0 or 1. All that is left to show is that the group 
\\

\setlength{\unitlength}{1.0cm}
\begin{picture}(12,1.0)
\thicklines
\put(1.4,0.1){$x_1$}
\put(3.0,0.1){$x_2$}
\put(4.6,0.1){$x_3$}
\put(6.2,0.1){$x_4$}
\put(7.8,0.1){$x_5$}
\put(9.4,0.1){$x_6$}
\put(11.0,0.1){$x_7$}
\put(2.2,0.6){$z_1$}
\put(3.8,0.6){$z_2$}
\put(5.4,0.6){$z_1$}
\put(7.0,0.6){$z_2$}
\put(10.0,0.6){$z_1 z_2$}
\put(1.9,0.2){\vector(1,0){1.0}}
\put(3.5,0.2){\vector(1,0){1.0}}
\put(5.1,0.2){\vector(1,0){1.0}}
\put(6.7,0.2){\vector(1,0){1.0}}
\put(9.9,0.2){\vector(1,0){1.0}}
\end{picture}
\\

\noindent
is isomorphic to A. Setting $x_3^\prime = x_3 x_1^{-1} x_5^{-1}$  and $z_1^\prime = z_1 z_2$ carries this group to
\\

\setlength{\unitlength}{1.0cm}
\begin{picture}(12,1.0)
\thicklines
\put(1.4,0.1){$x_1$}
\put(3.0,0.1){$x_2$}
\put(4.6,0.1){$x_3^\prime$}
\put(6.2,0.1){$x_4$}
\put(7.8,0.1){$x_5$}
\put(9.4,0.1){$x_6$}
\put(11.0,0.1){$x_7$}
\put(1.9,0.6){$z_1^\prime z_2^{-1}$}
\put(3.8,0.6){$z_1^\prime$}
\put(5.4,0.6){$z_1^\prime$}
\put(7.0,0.6){$z_2$}
\put(10.1,0.6){$z_1^\prime$}
\put(1.9,0.2){\vector(1,0){1.0}}
\put(3.5,0.2){\vector(1,0){1.0}}
\put(5.1,0.2){\vector(1,0){1.0}}
\put(6.7,0.2){\vector(1,0){1.0}}
\put(9.9,0.2){\vector(1,0){1.0}}
\end{picture}
\\

\noindent
Dropping primes and redrawing, this looks like the digraph below on the left. The substitution $x_4^\prime = x_4 x_2$ now carries this to the digraph below on the right.
\\

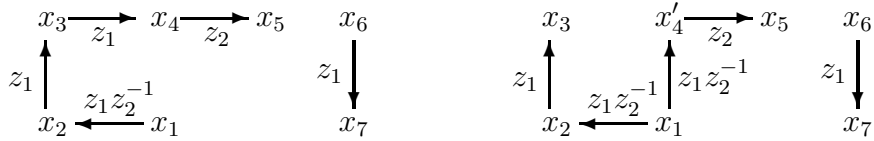
\setlength{\unitlength}{1.0cm}
\begin{picture}(12,2.0)
\thicklines
\put(2.3,0.1){$x_1$}      
\put(0.8,0.1){$x_2$}
\put(0.8,1.5){$x_3$}
\put(2.3,1.5){$x_4$}
\put(3.7,1.5){$x_5$}
\put(4.8,1.5){$x_6$}
\put(4.8,0.1){$x_7$}
\put(1.4,0.4){$z_1 z_2^{-1}$}
\put(0.4,0.7){$z_1$}
\put(1.5,1.3){$z_1$}
\put(3.0,1.3){$z_2$}
\put(4.5,0.8){$z_1$}
\put(2.2,0.2){\vector(-1,0){0.9}}
\put(0.9,0.4){\vector(0,1){0.9}}
\put(1.2,1.6){\vector(1,0){0.9}}
\put(2.7,1.6){\vector(1,0){0.9}}
\put(5.0,1.3){\vector(0,-1){0.9}}

\put(9.0,0.1){$x_1$}      
\put(7.5,0.1){$x_2$}
\put(7.5,1.5){$x_3$}
\put(9.0,1.5){$x_4^\prime$}
\put(10.4,1.5){$x_5$}
\put(11.5,1.5){$x_6$}
\put(11.5,0.1){$x_7$}
\put(8.1,0.4){$z_1 z_2^{-1}$}
\put(7.1,0.7){$z_1$}
\put(9.3,0.7){$z_1 z_2^{-1}$}
\put(9.7,1.3){$z_2$}
\put(11.2,0.8){$z_1$}
\put(8.9,0.2){\vector(-1,0){0.9}}
\put(7.6,0.4){\vector(0,1){0.9}}
\put(9.2,0.4){\vector(0,1){0.9}}
\put(9.4,1.6){\vector(1,0){0.9}}
\put(11.7,1.3){\vector(0,-1){0.9}}
\end{picture}
\\

\noindent
Setting $x_1^\prime = x_1 x_3 x_5^{-1}$ and inverting $x_3$, carries the group to 
\\

\setlength{\unitlength}{1.0cm}
\begin{picture}(12,1.0)
\thicklines
\put(1.4,0.1){$x_3$}
\put(3.1,0.1){$x_2$}
\put(4.8,0.1){$x_1$}
\put(6.5,0.1){$x_4$}
\put(8.2,0.1){$x_5$}
\put(9.6,0.1){$x_6$}
\put(11.3,0.1){$x_7$}
\put(2.2,0.6){$z_1$}
\put(4.0,0.6){$z_2$}
\put(5.8,0.6){$z_1$}
\put(7.6,0.6){$z_2$}
\put(10.3,0.6){$z_1$}
\put(1.9,0.2){\vector(1,0){1.0}}
\put(3.7,0.2){\vector(1,0){1.0}}
\put(5.5,0.2){\vector(1,0){1.0}}
\put(7.1,0.2){\vector(1,0){1.0}}
\put(10.0,0.2){\vector(1,0){1.0}}
\end{picture}
\\
\\
\noindent
Interchanging $x_1$  and $x_3$, this becomes group A. This completes the proof that there are six special groups of exponent $p$ having class 2 and derived group of order $p^2$.
\\

\section*{Acknowledgment}                                    
\noindent
The author would like to thank Michael Vaughan-Lee for his help in writing this paper and for his reading of Vishnevetskii's papers. Any errors or misstatements remaining in this paper are solely the responsibility of this author.
\\

\end{document}